\documentclass[10pt]{article}
\usepackage{cite}
\usepackage{mathrsfs}
\usepackage{amsfonts}
\usepackage{amsmath}
\usepackage{amsfonts,amssymb}
\usepackage{dsfont}
\usepackage{curves}
\usepackage{mathrsfs}
\usepackage{pifont}
\usepackage{amssymb}
\allowdisplaybreaks

\numberwithin{equation}{section}

\date{}

\def\BigRoman{\uppercase\expandafter{\romannumeral\number\count 255 }}
\def\Romannumeral{\afterassignment\BigRoman\count255=}

\begin{document}
\title{$A_{\alpha}$-spectral radius and path-factor covered graphs
\thanks{This work was supported by the Natural Science Foundation of Shandong Province, China (ZR2023MA078)}
}
\author{\small  Sizhong Zhou$^{1}$\footnote{Corresponding
author. E-mail address: zsz\_cumt@163.com (S. Zhou)}, Hongxia Liu$^{2}$, Qiuxiang Bian$^{1}$\\
\small $1$. School of Science, Jiangsu University of Science and Technology,\\
\small Zhenjiang, Jiangsu 212100, China\\
\small $2$. School of Mathematics and Information Science, Yantai University,\\
\small Yantai, Shandong 264005, China\\
}

\maketitle
\begin{abstract}
\noindent Let $\alpha\in[0,1)$, and let $G$ be a connected graph of order $n$ with $n\geq f(\alpha)$, where $f(\alpha)=14$ for
$\alpha\in[0,\frac{1}{2}]$, $f(\alpha)=17$ for $\alpha\in(\frac{1}{2},\frac{2}{3}]$, $f(\alpha)=20$ for $\alpha\in(\frac{2}{3},\frac{3}{4}]$
and $f(\alpha)=\frac{5}{1-\alpha}+1$ for $\alpha\in(\frac{3}{4},1)$. A path factor is a spanning subgraph $F$ of $G$ such that every component
of $F$ is a path with at least two vertices. Let $k\geq2$ be an integer. A $P_{\geq k}$-factor means a path-factor with each component being a
path of order at least $k$. A graph $G$ is called a $P_{\geq k}$-factor covered graph if $G$ has a $P_{\geq k}$-factor containing $e$ for any
$e\in E(G)$. Let $A_{\alpha}(G)=\alpha D(G)+(1-\alpha)A(G)$, where $D(G)$ denotes the diagonal matrix of vertex degrees of $G$ and $A(G)$ denotes
the adjacency matrix of $G$. The largest eigenvalue of $A_{\alpha}(G)$ is called the $A_{\alpha}$-spectral radius of $G$, which is denoted by
$\rho_{\alpha}(G)$. In this paper, it is proved that $G$ is a $P_{\geq2}$-factor covered graph if $\rho_{\alpha}(G)>\eta(n)$, where $\eta(n)$
is the largest root of $x^{3}-((\alpha+1)n+\alpha-4)x^{2}+(\alpha n^{2}+(\alpha^{2}-2\alpha-1)n-2\alpha+1)x-\alpha^{2}n^{2}+(5\alpha^{2}-3\alpha+2)n
-10\alpha^{2}+15\alpha-8=0$. Furthermore, we provide a graph to show that the bound on $A_{\alpha}$-spectral radius is optimal.
\\
\begin{flushleft}

(2020) Mathematics Subject Classification: 05C50, 05C70, 05C38
\end{flushleft}
\end{abstract}

\section{Introduction}

In this paper, we only consider finite and undirected graphs which have neither loops nor multiple edges. Let $G=(V(G),E(G))$ be a graph,
where $V(G)$ is the vertex set of $G$ and $E(G)$ is the edge set of $G$. The order of $G$ is denoted by $|V(G)|=n$. A graph $G$ is called
trivial if $|V(G)|=1$. Let $i(G)$ denote the number of isolated vertices in $G$. For any $v\in V(G)$, the neighborhood $N_G(v)$ of $v$ in
$G$ is defined by $\{u\in V(G):uv\in E(G)\}$. The cardinality
of $N_G(v)$ is called the degree of $v$ in $G$ and is denoted by $d_G(v)$. For a vertex subset $S$ of $G$, the induced subgraph induced by
$S$, denoted by $G[S]$, is the subgraph of $G$ whose vertex set is $S$ and whose edge set consists of all edges of $G$ which have both ends
in $S$. We write $G-S$ for $G[V(G)\setminus S]$. A subset $S\subseteq V(G)$ is called an independent set if $G[S]$ has no edges. The
complement of a graph $G$ is a graph $\overline{G}$ with the same vertex set as $G$, in which any two distinct vertices are adjacent if and
only if they are nonadjacent in $G$. Let $P_n$ and $K_n$ denote the path and complete graph of order $n$, respectively. Let $G_1$ and $G_2$
be two disjoint graphs. The union of $G_1$ and $G_2$, denoted by $G_1\cup G_2$, is the graph with vertex set $V(G_1)\cup V(G_2)$ and edge
set $E(G_1)\cup E(G_2)$. The join of $G_1$ and $G_2$, denoted by $G_1\vee G_2$, is the graph with vertex set $V(G_1)\cup V(G_2)$ and edge
set $E(G_1)\cup E(G_2)\cup\{uv:u\in V(G_1),v\in V(G_2)\}$. Let $k\geq3$ be an integer. The sequential join of graphs $G_1,G_2,\cdots,G_k$,
denoted by $G_1\vee G_2\vee\cdots\vee G_k$, is the graph with vertex set $V(G_1)\cup V(G_2)\cup\cdots\cup V(G_k)$ and edge set
$E(G_1)\cup E(G_2)\cup\cdots\cup E(G_k)\cup\{e=x_ix_{i+1}:x_i\in V(G_i), x_{i+1}\in V(G_{i+1}), 1\leq i\leq k-1\}$.

Given a graph $G$ with vertex set $V(G)=\{v_1,v_2,\cdots,v_n\}$, the adjacency matrix of $G$ is defined as $A(G)=(a_{ij})$, where $a_{ij}=1$
if $v_iv_j\in E(G)$, and $a_{ij}=0$ otherwise. Let $D(G)$ denote the diagonal matrix of vertex degrees of $G$. The signless Laplacian matrix
$Q(G)$ of $G$ are defined by $Q(G)=D(G)+A(G)$. For any $\alpha\in[0,1)$, Nikiforov \cite{N} defined the $A_{\alpha}$-matrix of $G$ as
$$
A_{\alpha}(G)=\alpha D(G)+(1-\alpha)A(G).
$$
Obviously, $A_0(G)=A(G)$ and $A_{\frac{1}{2}}(G)=\frac{1}{2}Q(G)$. Thus, $A_{\alpha}(G)$ generalizes both the adjacency matrix and the signless
Laplacian matrix of $G$. Note that $A_{\alpha}(G)$ is a real symmetric nonnegative matrix. Hence, the eigenvalues of $A_{\alpha}(G)$ are real.
The largest eigenvalue of $A_{\alpha}(G)$ is called the $A_{\alpha}$-spectral radius of $G$, which is denoted by $\rho_{\alpha}(G)$. More
contents on $A_{\alpha}$-matrix can be found in \cite{NR,So,LLX,BFO}.

A subset $M\subseteq E(G)$ is called a matching if any two edges of $M$ have no common vertices. $M$ is called a perfect matching if each
vertex of $G$ is incident to exactly one edge of $M$. A path factor is a spanning subgraph $F$ of $G$ such that every component of $F$ is a
path with at least two vertices. A $P_{\geq k}$-factor means a path-factor with every component having order at least $k$, where $k\geq2$ is
an integer. In fact, a graph $G$ has a perfect matching if and only if $G$ has a $\{P_2\}$-factor. A graph $G$ is called a $P_{\geq k}$-factor
covered graph if $G$ has a $P_{\geq k}$-factor containing $e$ for any $e\in E(G)$.

O \cite{O} obtained an adjacency spectral radius condition for the existence of a perfect matching in a graph. Liu, Pan and Li \cite{LPL} established
a relationship between a perfect matching and a signless Laplacian spectral radius of a graph. Zhang and Lin \cite{ZLp} proved an upper bound for
the distance spectral radius to guarantee the existence of a perfect matching in a graph. Zhao, Huang and Wang \cite{ZHW} presented an
$A_{\alpha}$-spectral radius condition for a graph to contain a perfect matching. For the relationship between spectral radius and the existence
of certain spanning subgraphs in graphs, one may be referred to \cite{ZLt,ZZL,FLL} and the references cited in. Egawa and Furuya \cite{EF} got a
sufficient condition for a graph to possess a $\{P_2,P_5\}$-factor. Wang and Zhang \cite{WZd} gave a degree condition for the existence of a
$\{P_2,P_5\}$-factor in a graph. Las Vergnas \cite{LV} showed a characterization for a graph with a $P_{\geq2}$-factor. Kaneko \cite{K} provided
a necessary and sufficient condition for a graph to contain a $P_{\geq3}$-factor. Kano, Lu and Yu \cite{KLY} established a connection between the
number of isolated vertices and $P_{\geq3}$-factors in graphs. Zhou and Sun \cite{ZS1} claimed a binding number condition for a graph with
a $P_{\geq2}$-factor (resp. a $P_{\geq3}$-factor). Wang and Zhang \cite{WZhi} presented a minimum degree and independence number condition for a
graph to contain a $P_{\geq3}$-factor. Wang and Zhang \cite{WZi} studied the relationship between isolated toughness and a $P_{\geq2}$-factor (resp.
a $P_{\geq3}$-factor) in a graph. Zhou, Sun and Bian \cite{ZSB} posed an isolated toughness condition for a graph to have a $P_{\geq3}$-factor. Zhou
\cite{Zs,Zd}, Zhou, Sun and Yang \cite{ZSY} obtained some results on the existence of $P_{\geq3}$-factors in graphs with given properties. Zhou, Sun
and Liu \cite{ZSL}, Wu \cite{Wp}, Liu \cite{Ls} gave some sufficient conditions for graphs to possess $P_{\geq3}$-factors with given properties.
More results on graph factors, we refer the reader to \cite{WZs,WZo,WZh,Wa,ZZs,ZPX1,ZPX,ZSL1}.

Zhou, Zhang and Sun \cite{ZZS} provided an $A_{\alpha}$-spectral radius condition for a graph with a $P_{\geq2}$-factor, which is shown in the
following.

\medskip

\noindent{\textbf{Theorem 1.1}} (Zhou, Zhang and Sun \cite{ZZS}). Let $\alpha\in[0,1)$, and let $G$ be a connected graph of order $n$ with $n\geq f(\alpha)$,
where
\[
f(\alpha)=\left\{
\begin{array}{ll}
14,&if \ \alpha\in[0,\frac{1}{2}];\\
17,&if \ \alpha\in(\frac{1}{2},\frac{2}{3}];\\
20,&if \ \alpha\in(\frac{2}{3},\frac{3}{4}];\\
\frac{5}{1-\alpha}+1,&if \ \alpha\in(\frac{3}{4},1).\\
\end{array}
\right.
\]
If $\rho_{\alpha}(G)>\theta(n)$, then $G$ contains a $P_{\geq2}$-factor, where $\theta(n)$ is the largest root of
$x^{3}-((\alpha+1)n+\alpha-5)x^{2}+(\alpha n^{2}+(\alpha^{2}-3\alpha-1)n-2\alpha+1)x-\alpha^{2}n^{2}+(7\alpha^{2}-5\alpha+3)n-18\alpha^{2}+29\alpha-15=0$.

\medskip

In this paper, we extend Theorem 1.1 to $P_{\geq2}$-factor covered graphs, and obtain a result on the existence of a $P_{\geq2}$-factor covered
graph by using $A_{\alpha}$-spectral radius.

\medskip

\noindent{\textbf{Theorem 1.2.}} Let $\alpha\in[0,1)$, and let $G$ be a connected graph of order $n$ with $n\geq f(\alpha)$, where $f(\alpha)$ is
defined as in Theorem 1.1. If $\rho_{\alpha}(G)>\eta(n)$, then $G$ is a $P_{\geq2}$-factor covered graph, where $\eta(n)$ is the largest root of
$x^{3}-((\alpha+1)n+\alpha-4)x^{2}+(\alpha n^{2}+(\alpha^{2}-2\alpha-1)n-2\alpha+1)x-\alpha^{2}n^{2}+(5\alpha^{2}-3\alpha+2)n
-10\alpha^{2}+15\alpha-8=0$.

\section{Preliminaries}

In this section, we provide several necessary preliminary results, which will be used to verify Theorem 1.2. Las Vergnas \cite{LV} posed a necessary
and sufficient condition for graphs to possess $P_{\geq2}$-factors.

\medskip

\noindent{\textbf{Lemma 2.1}} (Las Vergnas \cite{LV}). A graph $G$ has a $P_{\geq2}$-factor if and only if
$$
i(G-S)\leq2|S|
$$
for any vertex subset $S$ of $G$.

\medskip

Zhang and Zhou \cite{ZZ} obtained a necessary and sufficient condition for the existence of a $P_{\geq2}$-factor covered graph.

\medskip

\noindent{\textbf{Lemma 2.2}} (Zhang and Zhou \cite{ZZ}). Let $G$ be a connected graph. Then $G$ is a $P_{\geq2}$-factor covered graph if and
only if the following three conditions hold for any $S\subseteq V(G)$:

(\romannumeral1) $i(G-S)\leq2|S|$;

(\romannumeral2) If $S$ is nonempty and there exists a connected component of $G-S$ with at least two vertices, then $i(G-S)\leq2|S|-1$;

(\romannumeral3) If $S$ is not an independent set, then $i(G-S)\leq2|S|-2$.

\medskip

\noindent{\textbf{Lemma 2.3}} (Nikiforov \cite{N}). For a complete graph $K_n$ of order $n$, we have
$$
\rho_{\alpha}(K_n)=n-1.
$$

\medskip

\noindent{\textbf{Lemma 2.4}} (Nikiforov \cite{N}). Let $G$ be a connected graph, and $H$ be a proper subgraph of $G$. Then
$$
\rho_{\alpha}(G)>\rho_{\alpha}(H).
$$

\medskip

Let $M$ be a real symmetric matrix whose rows and columns are indexed by $V=\{1,2,\cdots,n\}$. Assume that $M$, with respect to the partition
$\pi: V=V_1\cup V_2\cup\cdots\cup V_r$, can be written as
\begin{align*}
M=\left(
  \begin{array}{ccc}
    M_{11} & \cdots & M_{1,r}\\
    \vdots & \ddots & \vdots\\
    M_{r1} & \cdots & M_{r,r}\\
  \end{array}
\right),
\end{align*}
where $M_{ij}$ denotes the submatrix (block) of $M$ formed by rows in $V_i$ and columns in $V_j$. Let $q_{ij}$ denote the average row sum of
$M_{ij}$. Then matrix $M_{\pi}=(q_{ij})$ is called the quotient matrix of $M$. If the row sum of each block $M_{ij}$ is a constant, then the
partition is equitable.

\medskip

\noindent{\textbf{Lemma 2.5}} (\cite{YYSX}).  Let $M$ be a real symmetric matrix with an equitable partition $\pi$, and
let $M_{\pi}$ be the corresponding quotient matrix. Then every eigenvalue of $M_{\pi}$ is an eigenvalue of $M$. Furthermore, if $M$ is
nonnegative, then the largest eigenvalues of $M$ and $M_{\pi}$ are equal.

\medskip

\section{The proof of Theorem 1.2}

\noindent{\it Proof of Theorem 1.2.} Suppose to the contrary that $G$ is not $P_{\geq2}$-factor covered. Choose a connected graph $G$ of order
$n$ such that its $A_{\alpha}$-spectral radius is as large as possible.

Let $h(x)=x^{3}-((\alpha+1)n+\alpha-5)x^{2}+(\alpha n^{2}+(\alpha^{2}-3\alpha-1)n-2\alpha+1)x-\alpha^{2}n^{2}+(7\alpha^{2}-5\alpha+3)n
-18\alpha^{2}+29\alpha-15$. By the condition of Theorem 1.1, $\theta(n)$ is the largest root of $h(x)=0$. Let
$\varphi(x)=x^{3}-((\alpha+1)n+\alpha-4)x^{2}+(\alpha n^{2}+(\alpha^{2}-2\alpha-1)n-2\alpha+1)x-\alpha^{2}n^{2}+(5\alpha^{2}-3\alpha+2)n
-10\alpha^{2}+15\alpha-8$. According to the condition of Theorem 1.2, $\eta(n)$ is the largest root of $\varphi(x)=0$. We verify Theorem 1.2 by
the next two possible cases.

\noindent{\bf Case 1.} $G$ has no $P_{\geq2}$-factor.

In terms of Theorem 1.1, we have $\rho_{\alpha}(G)\leq\theta(n)$ for $n\geq f(\alpha)$, where
\[
f(\alpha)=\left\{
\begin{array}{ll}
14,&if \ \alpha\in[0,\frac{1}{2}];\\
17,&if \ \alpha\in(\frac{1}{2},\frac{2}{3}];\\
20,&if \ \alpha\in(\frac{2}{3},\frac{3}{4}];\\
\frac{5}{1-\alpha}+1,&if \ \alpha\in(\frac{3}{4},1).\\
\end{array}
\right.
\]
It suffices to show $\theta:=\theta(n)\leq\eta(n)$. Notice that $h(\theta)=0$. Plugging $\theta$ into $x$ of $\varphi(x)-h(x)$ gives
$$
\varphi(\theta)=\varphi(\theta)-h(\theta)=-\theta^{2}+\alpha n\theta-(2\alpha^{2}-2\alpha+1)n+8\alpha^{2}-14\alpha+7.
$$

\noindent{\bf Claim 1.} $\theta>n-4$.

\noindent{\it Proof.} Let $H=K_1\vee(K_{n-4}\cup3K_1)$. We consider the partition $V(H)=V(3K_1)\cup V(K_{n-4})\cup V(K_1)$. The corresponding
quotient matrix of $A_{\alpha}(H)$ equals
\begin{align*}
B_0=\left(
  \begin{array}{ccc}
    \alpha & 0 & 1-\alpha\\
    0 & n+\alpha-5 & 1-\alpha\\
    3(1-\alpha) & (1-\alpha)(n-4) & \alpha n-\alpha\\
  \end{array}
\right).
\end{align*}
By a direct computation, the characteristic polynomial of $B_0$ is
\begin{align*}
h(x)=&x^{3}-((\alpha+1)n+\alpha-5)x^{2}+(\alpha n^{2}+(\alpha^{2}-3\alpha-1)n-2\alpha+1)x\\
&-\alpha^{2}n^{2}+(7\alpha^{2}-5\alpha+3)n-18\alpha^{2}+29\alpha-15.
\end{align*}
Notice that the partition $V(H)=V(3K_1)\cup V(K_{n-4})\cup V(K_1)$ is equitable. According to Lemma 2.5, the largest
root $\theta$ of $h(x)=0$ equals $\rho_{\alpha}(H)$, that is, $\rho_{\alpha}(H)=\theta$.

Obviously, $K_{n-3}$ is a proper subgraph of $H=K_1\vee(K_{n-4}\cup3K_1)$. Using Lemmas 2.3 and 2.4, we conclude
$$
\theta=\rho_{\alpha}(H)=\rho_{\alpha}(K_1\vee(K_{n-4}\cup3K_1))>\rho_{\alpha}(K_{n-3})=n-4.
$$
Claim 1 is true. \hfill $\Box$

By virtue of Claim 1, we get
$$
\frac{\alpha n}{2}<n-4<\theta.
$$
Then
\begin{align*}
\varphi(\theta)=&-\theta^{2}+\alpha n\theta-(2\alpha^{2}-2\alpha+1)n+8\alpha^{2}-14\alpha+7\\
<&-(n-4)^{2}+\alpha n(n-4)-(2\alpha^{2}-2\alpha+1)n+8\alpha^{2}-14\alpha+7\\
=&-(1-\alpha)n^{2}+(7-2\alpha-2\alpha^{2})n+8\alpha^{2}-14\alpha-9.
\end{align*}

If $0\leq\alpha\leq\frac{3}{4}$, then $\frac{7-2\alpha-2\alpha^{2}}{2-2\alpha}<14\leq f(\alpha)\leq n$. Thus, we possess
\begin{align*}
\varphi(\theta)<&-(1-\alpha)n^{2}+(7-2\alpha-2\alpha^{2})n+8\alpha^{2}-14\alpha-9\\
\leq&-196(1-\alpha)+14(7-2\alpha-2\alpha^{2})+8\alpha^{2}-14\alpha-9\\
=&-20\alpha^{2}+154\alpha-107\\
<&0.
\end{align*}

If $\frac{3}{4}<\alpha<1$, then $\frac{7-2\alpha-2\alpha^{2}}{2-2\alpha}<\frac{5}{1-\alpha}+1=f(\alpha)\leq n$. Thus, we obtain
\begin{align*}
\varphi(\theta)<&-(1-\alpha)n^{2}+(7-2\alpha-2\alpha^{2})n+8\alpha^{2}-14\alpha-9\\
\leq&-(1-\alpha)\left(\frac{5}{1-\alpha}+1\right)^{2}+(7-2\alpha-2\alpha^{2})\left(\frac{5}{1-\alpha}+1\right)+8\alpha^{2}-14\alpha-9\\
=&\frac{1}{1-\alpha}(-6\alpha^{3}+11\alpha^{2}-12\alpha-3)\\
<&0.
\end{align*}

From the discussion above, we always have $\varphi(\theta)<0$, which implies $\theta(n)=\theta<\eta(n)$. Together with
$\rho_{\alpha}(G)\leq\theta(n)$ for $n\geq f(\alpha)$, we conclude $\rho_{\alpha}(G)\leq\theta(n)<\eta(n)$  for $n\geq f(\alpha)$, which is a
contradiction to $\rho_{\alpha}(G)>\eta(n)$.

\noindent{\bf Case 2.} $G$ has a $P_{\geq2}$-factor.

According to Lemma 2.1, we obtain $i(G-S)\leq2|S|$ for any subset $S$ of $V(G)$. Since $G$ is not $P_{\geq2}$-factor covered, using Lemma 2.2,
there exists a subset $S$ of $V(G)$ so that either $S$ is nonempty and some connected component of $G-S$ is nontrivial, or $S$ is not an independent
set. Otherwise, it follows from Lemma 2.2 that $G$ is $P_{\geq2}$-factor covered, a contradiction. Furthermore, at least one of the following
statements holds:

\noindent{(a)} $S$ is nonempty, one connected component of $G-S$ is nontrivial, and $i(G-S)\geq2|S|$;

\noindent{(b)} $S$ is not an independent set, and $i(G-S)\geq2|S|-1$.

\noindent{\bf Subcase 2.1.} Statement (a) holds.

In this subcase, we have $|S|\geq1$ and $i(G-S)=2|S|$. In view of Lemma 2.4 and the choice of $G$, we conclude that

$\bullet$ the induced subgraph $G[S]$ and every connected component of $G-S$ are complete graphs;

$\bullet$ $G=G[S]\vee(G-S)$;

$\bullet$ $G-S$ has only one nontrivial connected component, say $G_1$.

Let $|S|=s$ and $|V(G_1)|=n_1\geq2$. Recall that $i(G-S)=2|S|$. Then $G=K_s\vee(K_{n_1}\cup2sK_1)$, where $n_1=n-3s\geq2$. It is obvious that
the partition $V(G)=V(K_s)\cup V(K_{n_1})\cup V(2sK_1)$ is an equitable partition of $G$, and the corresponding quotient matrix of $A_{\alpha}(G)$
is
\begin{align*}
B_1=\left(
  \begin{array}{ccc}
    \alpha n-\alpha s+s-1 & (1-\alpha)(n-3s) & 2s(1-\alpha)\\
    (1-\alpha)s & n+\alpha s-3s-1 & 0\\
    (1-\alpha)s & 0 & \alpha s\\
  \end{array}
\right).
\end{align*}
By a simple computation, the characteristic polynomial of $B_1$ is given by
\begin{align}\label{eq:3.1}
\varphi_1(x)=&x^{3}-((\alpha+1)n+(\alpha-2)s-2)x^{2}\nonumber\\
&+(\alpha n^{2}+(\alpha^{2}-\alpha)sn-(\alpha+1)n-2s^{2}-(2\alpha-2)s+1)x\nonumber\\
&-\alpha^{2}sn^{2}+(4\alpha^{2}-4\alpha+2)s^{2}n+(\alpha^{2}+\alpha)sn\nonumber\\
&-(8\alpha^{2}-14\alpha+6)s^{3}-(2\alpha^{2}-2\alpha+2)s^{2}-\alpha s.
\end{align}
By virtue of Lemma 2.5, $\rho_{\alpha}(G)$ is the largest root of $\varphi_1(x)=0$. Let $\eta_1=\rho_{\alpha}(G)\geq\eta_2\geq\eta_3$ be the three
roots of $\varphi_1(x)=0$ and $R_1=diag(s,n-3s,2s)$. One checks that $R_1^{\frac{1}{2}}B_1R_1^{-\frac{1}{2}}$ is symmetric, and also contains
\begin{align*}
\left(
  \begin{array}{ccc}
    n+\alpha s-3s-1 & 0\\
    0 &  \alpha s\\
  \end{array}
\right)
\end{align*}
as its submatrix. Since $R_1^{\frac{1}{2}}B_1R_1^{-\frac{1}{2}}$ and $B_1$ possess the same eigenvalues, the Cauchy interlacing theorem (cf. \cite{H})
implies that $\eta_2\leq n+\alpha s-3s-1<n-3$.

If $s=1$, then $G=K_1\vee(K_{n-3}\cup2K_1)$ and $\varphi_1(x)=\varphi(x)$. Recall that $\eta(n)$ is the largest root of $\varphi(x)=0$ and
$\rho_{\alpha}(G)$ is the largest root of $\varphi_1(x)=0$. Thus, we conclude $\rho_{\alpha}(G)=\eta(n)$, which is a contradiction to
$\rho_{\alpha}(G)>\eta(n)$.  In what follows, we are to consider $s\geq2$.

From the discussion above, we easily see $\eta(n)=\rho_{\alpha}(K_1\vee(K_{n-3}\cup2K_1))$. Note that $K_{n-2}$ is a proper subgraph of
$K_1\vee(K_{n-3}\cup2K_1)$. From Lemmas 2.3 and 2.4, we have $\rho_{\alpha}(K_1\vee(K_{n-3}\cup2K_1))>\rho_{\alpha}(K_{n-2})=n-3$. Thus, we obtain
\begin{align}\label{eq:3.2}
\eta(n)=\rho_{\alpha}(K_1\vee(K_{n-3}\cup2K_1))>\rho_{\alpha}(K_{n-2})=n-3>\eta_2.
\end{align}

\noindent{\bf Subcase 2.1.1.} $0\leq\alpha\leq\frac{3}{4}$.

Write $\eta:=\eta(n)$. Recall that $\varphi(\eta)=0$. A directed calculation yields that
\begin{align}\label{eq:3.3}
\varphi_1(\eta)=\varphi_1(\eta)-\varphi(\eta)=(s-1)g_1(\eta),
\end{align}
where $g_1(\eta)=(2-\alpha)\eta^{2}-((\alpha-\alpha^{2})n+2s+2\alpha)\eta-\alpha^{2}n^{2}+(4\alpha^{2}-4\alpha+2)sn+(5\alpha^{2}-3\alpha+2)n
-(8\alpha^{2}-14\alpha+6)s^{2}-(10\alpha^{2}-16\alpha+8)s-10\alpha^{2}+15\alpha-8$. Note that
$$
\frac{(\alpha-\alpha^{2})n+2s+2\alpha}{2(2-\alpha)}<n-3<\eta
$$
by \eqref{eq:3.2}, $s\geq2$ and $n\geq3s+2$. Thus, we obtain
\begin{align}\label{eq:3.4}
g_1(\eta)>&g_1(n-3)\nonumber\\
=&(2-\alpha)(n-3)^{2}-((\alpha-\alpha^{2})n+2s+2\alpha)(n-3)-\alpha^{2}n^{2}\nonumber\\
&+(4\alpha^{2}-4\alpha+2)sn+(5\alpha^{2}-3\alpha+2)n-(8\alpha^{2}-14\alpha+6)s^{2}\nonumber\\
&-(10\alpha^{2}-16\alpha+8)s-10\alpha^{2}+15\alpha-8\nonumber\\
=&(2-2\alpha)n^{2}+((4\alpha^{2}-4\alpha)s+2\alpha^{2}+4\alpha-10)n-(8\alpha^{2}-14\alpha+6)s^{2}\nonumber\\
&-(10\alpha^{2}-16\alpha+2)s-10\alpha^{2}+12\alpha+10.
\end{align}
Let $l_1(s,n)=(2-2\alpha)n^{2}+((4\alpha^{2}-4\alpha)s+2\alpha^{2}+4\alpha-10)n-(8\alpha^{2}-14\alpha+6)s^{2}-(10\alpha^{2}-16\alpha+2)s-10\alpha^{2}+12\alpha+10$.
Notice that
$$
-\frac{(4\alpha^{2}-4\alpha)s+2\alpha^{2}+4\alpha-10}{2(2-2\alpha)}<3s+2\leq n.
$$
Then
\begin{align*}
l_1(s,n)\geq&l_1(s,3s+2)\\
=&(2-2\alpha)(3s+2)^{2}+((4\alpha^{2}-4\alpha)s+2\alpha^{2}+4\alpha-10)(3s+2)\\
&-(8\alpha^{2}-14\alpha+6)s^{2}-(10\alpha^{2}-16\alpha+2)s-10\alpha^{2}+12\alpha+10\\
=&(4s^{2}+4s-6)\alpha^{2}-(16s^{2}+4s-12)\alpha+12s^{2}-8s-2\\
\geq&\frac{9}{16}(4s^{2}+4s-6)-\frac{3}{4}(16s^{2}+4s-12)+12s^{2}-8s-2\\
=&\frac{1}{8}(18s^{2}-70s+29)\\
>&0,
\end{align*}
where the last two inequalities hold from $\frac{16s^{2}+4s-12}{2(4s^{2}+4s-6)}>\frac{3}{4}\geq\alpha$ and $s\geq4$, respectively.

If $s=3$, then $-\frac{(4\alpha^{2}-4\alpha)s+2\alpha^{2}+4\alpha-10}{2(2-2\alpha)}=-\frac{7\alpha^{2}-4\alpha-5}{2-2\alpha}<14\leq f(\alpha)\leq n$
due to $0\leq\alpha\leq\frac{3}{4}$. Consequently, for $0\leq\alpha\leq\frac{3}{4}$ and $n\geq f(\alpha)\geq14$, we obtain
$$
l_1(3,n)\geq l_1(3,14)=84\alpha^{2}-318\alpha+202>0.
$$

If $s=2$, then $-\frac{(4\alpha^{2}-4\alpha)s+2\alpha^{2}+4\alpha-10}{2(2-2\alpha)}=-\frac{5\alpha^{2}-2\alpha-5}{2-2\alpha}<14\leq f(\alpha)\leq n$
due to $0\leq\alpha\leq\frac{3}{4}$. Hence, for $0\leq\alpha\leq\frac{3}{4}$ and $n\geq f(\alpha)\geq14$, we conclude
$$
l_1(2,n)\geq l_1(2,14)=78\alpha^{2}-328\alpha+234>0.
$$

From the discussion above, we have $l_1(s,n)>0$ for $2\leq s\leq\frac{n-2}{3}$. Together with \eqref{eq:3.3}, \eqref{eq:3.4} and $s\geq2$, we get
$$
\varphi_1(\eta)=(s-1)g_1(\eta)>(s-1)g_1(n-3)=(s-1)l_1(s,n)>0.
$$
Recall that $\eta=\eta(n)$ and $\rho_{\alpha}(G)$ is the largest root of $\varphi_1(x)=0$. As $\eta(n)=\rho_{\alpha}(K_1\vee(K_{n-3}\cup2K_1))>n-3>\eta_2$
(see \eqref{eq:3.2}), we infer $\rho_{\alpha}(G)<\eta(n)$ for $2\leq s\leq\frac{n-2}{3}$, which is a contradiction to $\rho_{\alpha}(G)>\eta(n)$.

\noindent{\bf Subcase 2.1.2.} $\frac{3}{4}<\alpha<1$.

According to \eqref{eq:3.1}, we have
\begin{align*}
\varphi_1(n-3)=&(n-3)^{3}-((\alpha+1)n+(\alpha-2)s-2)(n-3)^{2}\\
&+(\alpha n^{2}+(\alpha^{2}-\alpha)sn-(\alpha+1)n-2s^{2}-(2\alpha-2)s+1)(n-3)\\
&-\alpha^{2}sn^{2}+(4\alpha^{2}-4\alpha+2)s^{2}n+(\alpha^{2}+\alpha)sn\\
&-(8\alpha^{2}-14\alpha+6)s^{3}-(2\alpha^{2}-2\alpha+2)s^{2}-\alpha s\\
=&2(1-\alpha)(4\alpha-3)s^{3}+((4\alpha^{2}-4\alpha)n-2\alpha^{2}+2\alpha+4)s^{2}\\
&+((2-2\alpha)n^{2}-(2\alpha^{2}-8\alpha+10)n-4\alpha+12)s\\
&+(2\alpha-2)n^{2}-(6\alpha-10)n-12\\
=:&\Phi_1(s,n).
\end{align*}
Thus, we obtain
\begin{align*}
\frac{\partial \Phi_1(s,n)}{\partial s}=&6(1-\alpha)(4\alpha-3)s^{2}+2((4\alpha^{2}-4\alpha)n-2\alpha^{2}+2\alpha+4)s\\
&+(2-2\alpha)n^{2}-(2\alpha^{2}-8\alpha+10)n-4\alpha+12.
\end{align*}
Recall that $\frac{3}{4}<\alpha<1$ and $n\geq f(\alpha)=\frac{5}{1-\alpha}+1>\frac{5}{1-\alpha}$. By a direct calculation, we have
\begin{align*}
\frac{\partial \Phi_1(s,n)}{\partial s}\Big|_{s=2}=&(2-2\alpha)n^{2}+(14\alpha^{2}-8\alpha-10)n-104\alpha^{2}+172\alpha-44\\
>&(2-2\alpha)\left(\frac{5}{1-\alpha}\right)^{2}+(14\alpha^{2}-8\alpha-10)\left(\frac{5}{1-\alpha}\right)-104\alpha^{2}+172\alpha-44\\
=&\frac{1}{1-\alpha}(104\alpha^{3}-206\alpha^{2}+176\alpha-44)\\
>&0,
\end{align*}
and
\begin{align*}
\frac{\partial \Phi_1(s,n)}{\partial s}\Big|_{s=\frac{n-2}{3}}=&\frac{1}{3}((6\alpha^{2}-12\alpha+2)n-24\alpha^{2}+36\alpha-4)\\
<&\frac{1}{3}\left((6\alpha^{2}-12\alpha+2)\left(\frac{5}{1-\alpha}\right)-24\alpha^{2}+36\alpha-4\right)\\
=&\frac{1}{3(1-\alpha)}(24\alpha^{3}-30\alpha^{2}-20\alpha+6)\\
<&0.
\end{align*}
This implies that $\varphi_1(n-3)=\Phi_1(s,n)\geq\min\left\{\Phi_1(2,n),\Phi_1\left(\frac{n-2}{3},n\right)\right\}$ because the leading coefficient
of $\Phi_1(s,n)$ (view as a cubic polynomial of $s$) is positive, and $2\leq s\leq\frac{n-2}{3}$. By virtue of $\frac{3}{4}<\alpha<1$ and
$n\geq f(\alpha)=\frac{5}{1-\alpha}+1>\frac{5}{1-\alpha}$, we obtain
\begin{align*}
\Phi_1(2,n)=&(2-2\alpha)n^{2}+(12\alpha^{2}-6\alpha-10)n-72\alpha^{2}+112\alpha-20\\
>&(2-2\alpha)\left(\frac{5}{1-\alpha}\right)^{2}+(12\alpha^{2}-6\alpha-10)\left(\frac{5}{1-\alpha}\right)-72\alpha^{2}+112\alpha-20\\
=&\frac{1}{1-\alpha}(72\alpha^{3}-124\alpha^{2}+102\alpha-20)\\
>&0,
\end{align*}
and
\begin{align*}
\Phi_1\left(\frac{n-2}{3},n\right)=&\frac{1}{27}((4\alpha^{2}-16\alpha+12)n^{3}+(-24\alpha^{2}+132\alpha-132)n^{2}\\
&+(12\alpha^{2}-246\alpha+438)n+40\alpha^{2}-16\alpha-444)\\
>&\frac{1}{27}((4\alpha^{2}-16\alpha+12)\left(\frac{5}{1-\alpha}\right)^{3}+(-24\alpha^{2}+132\alpha-132)\left(\frac{5}{1-\alpha}\right)^{2}\\
&+(12\alpha^{2}-246\alpha+438)\left(\frac{5}{1-\alpha}\right)+40\alpha^{2}-16\alpha-444)\\
=&\frac{1}{27(1-\alpha)^{2}}(40\alpha^{4}-156\alpha^{3}+318\alpha^{2}+252\alpha-54)\\
>&0.
\end{align*}
Consequently, we conclude $\varphi_1(n-3)\geq\min\left\{\Phi_1(2,n),\Phi_1\left(\frac{n-2}{3},n\right)\right\}>0$ for $2\leq s\leq\frac{n-2}{3}$. As
$\eta(n)=\rho_{\alpha}(K_1\vee(K_{n-3}\cup2K_1))>n-3>\eta_2$ (see \eqref{eq:3.2}), we deduce $\rho_{\alpha}(G)<\eta(n)$ for $2\leq s\leq\frac{n-2}{3}$,
which is a contradiction to $\rho_{\alpha}(G)>\eta(n)$.

\noindent{\bf Subcase 2.2.} Statement (b) holds.

In this subcase, we infer $|S|\geq2$ and $i(G-S)=2|S|-1 \ \mbox{or} \ 2|S|$. Similarly, one has

$\bullet$ the induced subgraph $G[S]$ and every connected component of $G-S$ are complete graphs;

$\bullet$ $G=G[S]\vee(G-S)$;

$\bullet$ $G-S$ has at most one nontrivial connected component $G_1$.

We also let $|S|=s\geq2$ and $|V(G_1)|=n_1$. Notice that $n_1=0$ or $n_1\geq2$.

\noindent{\bf Subcase 2.2.1.} $i(G-S)=2s-1$.

\noindent{\bf Subcase 2.2.1.1.} $n_1\geq2$.

In this subcase, $G=K_s\vee(K_{n_1}\cup(2s-1)K_1)$ and $n=3s-1+n_1\geq3s+1$. It is clear that the partition $V(G)=V(K_s)\cup V(K_{n_1})\cup V((2s-1)K_1)$
is an equitable partition of $G$, and the corresponding quotient matrix of $A_{\alpha}(G)$ equals
\begin{align*}
B_2=\left(
  \begin{array}{ccc}
    \alpha n-\alpha s+s-1 & (1-\alpha)(n-3s+1) & (1-\alpha)(2s-1)\\
    (1-\alpha)s & n+\alpha s-3s & 0\\
    (1-\alpha)s & 0 & \alpha s\\
  \end{array}
\right).
\end{align*}
By a direct computation, we obtain the characteristic polynomial of $B_2$ as
\begin{align}\label{eq:3.5}
\varphi_2(x)=&x^{3}-((\alpha+1)n+(\alpha-2)s-1)x^{2}\nonumber\\
&+(\alpha n^{2}+(\alpha^{2}-\alpha)sn-n-2s^{2}-(2\alpha-3)s)x\nonumber\\
&-\alpha^{2}sn^{2}+(4\alpha^{2}-4\alpha+2)s^{2}n-(\alpha^{2}-3\alpha+1)sn\nonumber\\
&-(8\alpha^{2}-14\alpha+6)s^{3}+(4\alpha^{2}-9\alpha+3)s^{2}.
\end{align}
According to Lemma 2.5, $\rho_{\alpha}(G)$ is the largest root of $\varphi_2(x)=0$. Let $\beta_1=\rho_{\alpha}(G)\geq\beta_2\geq\beta_3$ be the three
roots of $\varphi_2(x)=0$ and $R_2=diag(s,n-3s+1,2s-1)$. One checks that $R_2^{\frac{1}{2}}B_2R_2^{-\frac{1}{2}}$ is symmetric, and also contains
\begin{align*}
\left(
  \begin{array}{ccc}
    n+\alpha s-3s & 0\\
    0 &  \alpha s\\
  \end{array}
\right)
\end{align*}
as its submatrix. Since $R_2^{\frac{1}{2}}B_2R_2^{-\frac{1}{2}}$ and $B_2$ have the same eigenvalues, the Cauchy interlacing theorem (cf. \cite{H})
yields that $\beta_2\leq n+\alpha s-3s<n-2s\leq n-4$.

Recall that $\eta=\eta(n)=\rho_{\alpha}(K_1\vee(K_{n-3}\cup2K_1))$. Since $K_{n-2}$ is a proper subgraph of $K_1\vee(K_{n-3}\cup2K_1)$, it follows
from Lemmas 2.3 and 2.4 that
\begin{align}\label{eq:3.6}
\eta=\eta(n)=\rho_{\alpha}(K_1\vee(K_{n-3}\cup2K_1))>\rho_{\alpha}(K_{n-2})=n-3>n-4>\beta_2.
\end{align}

First we discuss $0\leq\alpha\leq\frac{3}{4}$. Recall that $\varphi(\eta)=0$. A simple computation yields that
\begin{align*}
\varphi_2(\eta)=&\varphi_2(\eta)-\varphi(\eta)\\
=&((2-\alpha)s+\alpha-3)\eta^{2}+(((\alpha^{2}-\alpha)s-\alpha^{2}+2\alpha)n-2s^{2}-(2\alpha-3)s+2\alpha-1)\eta\\
&-\alpha^{2}(s-1)n^{2}+(4\alpha^{2}-4\alpha+2)s^{2}n-(\alpha^{2}-3\alpha+1)sn-(5\alpha^{2}-3\alpha+2)n\\
&-(8\alpha^{2}-14\alpha+6)s^{3}+(4\alpha^{2}-9\alpha+3)s^{2}+10\alpha^{2}-15\alpha+8.
\end{align*}
Note that
$$
-\frac{((\alpha^{2}-\alpha)s-\alpha^{2}+2\alpha)n-2s^{2}-(2\alpha-3)s+2\alpha-1}{2((2-\alpha)s+\alpha-3)}<n-3<\eta
$$
by \eqref{eq:3.6}, $s\geq2$ and $n\geq3s+1$. Hence, we obtain
\begin{align}\label{eq:3.7}
\varphi_2(\eta)>&((2-\alpha)s+\alpha-3)(n-3)^{2}+(((\alpha^{2}-\alpha)s-\alpha^{2}+2\alpha)n-2s^{2}-(2\alpha-3)s+2\alpha-1)(n-3)\nonumber\\
&-\alpha^{2}(s-1)n^{2}+(4\alpha^{2}-4\alpha+2)s^{2}n-(\alpha^{2}-3\alpha+1)sn-(5\alpha^{2}-3\alpha+2)n\nonumber\\
&-(8\alpha^{2}-14\alpha+6)s^{3}+(4\alpha^{2}-9\alpha+3)s^{2}+10\alpha^{2}-15\alpha+8\nonumber\\
=&((2-2\alpha)s+3\alpha-3)n^{2}+((4\alpha^{2}-4\alpha)s^{2}-(4\alpha^{2}-10\alpha+10)s-2\alpha^{2}-7\alpha+15)n\nonumber\\
&-(8\alpha^{2}-14\alpha+6)s^{3}+(4\alpha^{2}-9\alpha+9)s^{2}-3(\alpha-3)s+10\alpha^{2}-12\alpha-16.
\end{align}
Let $l_2(s,n)=((2-2\alpha)s+3\alpha-3)n^{2}+((4\alpha^{2}-4\alpha)s^{2}-(4\alpha^{2}-10\alpha+10)s-2\alpha^{2}-7\alpha+15)n -(8\alpha^{2}-14\alpha+6)s^{3}+(4\alpha^{2}-9\alpha+9)s^{2}-3(\alpha-3)s+10\alpha^{2}-12\alpha-16$. Note that $(2-2\alpha)s+3\alpha-3>0$ and
$$
-\frac{(4\alpha^{2}-4\alpha)s^{2}-(4\alpha^{2}-10\alpha+10)s-2\alpha^{2}-7\alpha+15}{2((2-2\alpha)s+3\alpha-3)}<3s+1\leq n.
$$
Then
\begin{align*}
l_2(s,n)\geq&l_2(s,3s+1)\\
=&((2-2\alpha)s+3\alpha-3)(3s+1)^{2}\\
&+((4\alpha^{2}-4\alpha)s^{2}-(4\alpha^{2}-10\alpha+10)s-2\alpha^{2}-7\alpha+15)(3s+1)\\
&-(8\alpha^{2}-14\alpha+6)s^{3}+(4\alpha^{2}-9\alpha+9)s^{2}-3(\alpha-3)s+10\alpha^{2}-12\alpha-16\\
=&(4s^{3}-4s^{2}-10s+8)\alpha^{2}-(16s^{3}-32s^{2}-2s+16)\alpha+12s^{3}-36s^{2}+28s-4\\
\geq&\frac{9}{16}(4s^{3}-4s^{2}-10s+8)-\frac{3}{4}(16s^{3}-32s^{2}-2s+16)+12s^{3}-36s^{2}+28s-4\\
=&\frac{1}{8}(18s^{3}-114s^{2}+191s-92)\\
\geq&0,
\end{align*}
where the last two inequalities hold from $\frac{16s^{2}+4s-12}{2(4s^{2}+4s-6)}>\frac{3}{4}\geq\alpha$ and $s\geq4$, respectively.

If $s=3$, then $-\frac{(4\alpha^{2}-4\alpha)s^{2}-(4\alpha^{2}-10\alpha+10)s-2\alpha^{2}-7\alpha+15}{2((2-2\alpha)s+3\alpha-3)} =-\frac{22\alpha^{2}-13\alpha-15}{6-6\alpha}<14\leq f(\alpha)\leq n$ due to $0\leq\alpha\leq\frac{3}{4}$. Hence, for $0\leq\alpha\leq\frac{3}{4}$ and $n\geq f(\alpha)\geq14$, we get
$$
l_2(3,n)\geq l_2(3,14)=138\alpha^{2}-494\alpha+308>0.
$$

If $s=2$, then $-\frac{(4\alpha^{2}-4\alpha)s^{2}-(4\alpha^{2}-10\alpha+10)s-2\alpha^{2}-7\alpha+15}{2((2-2\alpha)s+3\alpha-3)}= -\frac{6\alpha^{2}-3\alpha-5}{2-2\alpha}<14\leq f(\alpha)\leq n$ due to $0\leq\alpha\leq\frac{3}{4}$. Hence, for $0\leq\alpha\leq\frac{3}{4}$ and $n\geq f(\alpha)\geq14$, we have
$$
l_2(2,n)\geq l_2(2,14)=46\alpha^{2}-180\alpha+115>0.
$$

From the discussion above, we conclude $l_2(s,n)\geq0$ for $2\leq s\leq\frac{n-1}{3}$. Combining this with \eqref{eq:3.7} and $s\geq2$, we obtain
$$
\varphi_2(\eta)>l_2(s,n)\geq0.
$$
Recall that $\rho_{\alpha}(G)$ is the largest root of $\varphi_2(x)=0$. As $\eta=\eta(n)=\rho_{\alpha}(K_1\vee(K_{n-3}\cup2K_1))>n-3>\beta_2$
(see \eqref{eq:3.6}), we deduce $\rho_{\alpha}(G)<\eta(n)$ for $2\leq s\leq\frac{n-1}{3}$, which is a contradiction to $\rho_{\alpha}(G)>\eta(n)$. In
what follows, we consider $\frac{3}{4}<\alpha<1$.

In terms of \eqref{eq:3.5}, we et
\begin{align*}
\varphi_2(n-3)=&(n-3)^{3}-((\alpha+1)n+(\alpha-2)s-1)(n-3)^{2}\\
&+(\alpha n^{2}+(\alpha^{2}-\alpha)sn-n-2s^{2}-(2\alpha-3)s)(n-3)\\
&-\alpha^{2}sn^{2}+(4\alpha^{2}-4\alpha+2)s^{2}n-(\alpha^{2}-3\alpha+1)sn\\
&-(8\alpha^{2}-14\alpha+6)s^{3}+(4\alpha^{2}-9\alpha+3)s^{2}\\
=&2(1-\alpha)(4\alpha-3)s^{3}+((4\alpha^{2}-4\alpha)n+4\alpha^{2}-9\alpha+9)s^{2}\\
&+((2-2\alpha)n^{2}-(4\alpha^{2}-10\alpha+10)n-3\alpha+9)s\\
&+(3\alpha-3)n^{2}-(9\alpha-15)n-18\\
=:&\Phi_2(s,n).
\end{align*}
Thus, we have
\begin{align*}
\frac{\partial \Phi_2(s,n)}{\partial s}=&6(1-\alpha)(4\alpha-3)s^{2}+2((4\alpha^{2}-4\alpha)n+4\alpha^{2}-9\alpha+9)s\\
&+(2-2\alpha)n^{2}-(4\alpha^{2}-10\alpha+10)n-3\alpha+9.
\end{align*}
Recall that $\frac{3}{4}<\alpha<1$ and $n\geq f(\alpha)=\frac{5}{1-\alpha}+1>\frac{5}{1-\alpha}$. By a simple calculation, we obtain
\begin{align*}
\frac{\partial \Phi_2(s,n)}{\partial s}\Big|_{s=2}=&(2-2\alpha)n^{2}+(12\alpha^{2}-6\alpha-10)n-80\alpha^{2}+129\alpha-27\\
>&(2-2\alpha)\left(\frac{5}{1-\alpha}\right)^{2}+(12\alpha^{2}-6\alpha-10)\left(\frac{5}{1-\alpha}\right)-80\alpha^{2}+129\alpha-27\\
=&\frac{1}{1-\alpha}(80\alpha^{3}-149\alpha^{2}+126\alpha-27)\\
>&0,
\end{align*}
and
\begin{align*}
\frac{\partial \Phi_2(s,n)}{\partial s}\Big|_{s=\frac{n-1}{3}}=&\frac{1}{3}((4\alpha^{2}-8\alpha)n-16\alpha^{2}+23\alpha+3)\\
<&\frac{1}{3}\left((4\alpha^{2}-8\alpha)\left(\frac{5}{1-\alpha}\right)-16\alpha^{2}+23\alpha+3\right)\\
=&\frac{1}{3(1-\alpha)}(16\alpha^{3}-19\alpha^{2}-20\alpha+3)\\
<&0.
\end{align*}
This implies that $\varphi_2(n-3)=\Phi_2(s,n)\geq\min\left\{\Phi_2(2,n),\Phi_2\left(\frac{n-1}{3},n\right)\right\}$ because the leading coefficient
of $\Phi_2(s,n)$ (view as a cubic polynomial of $s$) is positive, and $2\leq s\leq\frac{n-1}{3}$. In light of $\frac{3}{4}<\alpha<1$ and
$n\geq f(\alpha)=\frac{5}{1-\alpha}+1>\frac{5}{1-\alpha}$, we have
\begin{align*}
\Phi_2(2,n)=&(1-\alpha)n^{2}+(8\alpha^{2}-4\alpha-5)n-48\alpha^{2}+70\alpha-12\\
>&(1-\alpha)\left(\frac{5}{1-\alpha}\right)^{2}+(8\alpha^{2}-4\alpha-5)\left(\frac{5}{1-\alpha}\right)-48\alpha^{2}+70\alpha-12\\
=&\frac{1}{1-\alpha}(48\alpha^{3}-78\alpha^{2}+62\alpha-12)\\
>&0,
\end{align*}
and
\begin{align*}
\Phi_2\left(\frac{n-1}{3},n\right)=&\frac{1}{27}((4\alpha^{2}-16\alpha+12)n^{3}+(-24\alpha^{2}+144\alpha-144)n^{2}\\
&+(-249\alpha+504)n+20\alpha^{2}-14\alpha-534)\\
>&\frac{1}{27}((4\alpha^{2}-16\alpha+12)\left(\frac{5}{1-\alpha}\right)^{3}+(-24\alpha^{2}+144\alpha-144)\left(\frac{5}{1-\alpha}\right)^{2}\\
&+(-249\alpha+504)\left(\frac{5}{1-\alpha}\right)+20\alpha^{2}-14\alpha-534)\\
=&\frac{1}{27(1-\alpha)^{2}}(20\alpha^{4}-54\alpha^{3}+159\alpha^{2}+389\alpha-114)\\
>&0.
\end{align*}
Hence, we have $\varphi_2(n-3)\geq\min\left\{\Phi_2(2,n),\Phi_2\left(\frac{n-1}{3},n\right)\right\}>0$ for $2\leq s\leq\frac{n-1}{3}$. As
$\eta=\eta(n)=\rho_{\alpha}(K_1\vee(K_{n-3}\cup2K_1))>n-3>\beta_2$ (see \eqref{eq:3.6}), we conclude $\rho_{\alpha}(G)<\eta(n)$ for $2\leq s\leq\frac{n-1}{3}$,
which is a contradiction to $\rho_{\alpha}(G)>\eta(n)$.

\noindent{\bf Subcase 2.2.1.2.} $n_1=0$.

In this subcase, $G=K_s\vee(2s-1)K_1$ and $n=3s-1$. The quotient matrix of $A_{\alpha}(G)$ by the partition $V(G)=V(K_s)\cup V((2s-1)K_1)$ equals
\begin{align*}
B_3=\left(
  \begin{array}{ccc}
    2\alpha s+s-\alpha-1 & (1-\alpha)(2s-1)\\
    (1-\alpha)s &  \alpha s\\
  \end{array}
\right).
\end{align*}

By a direct computation, we get the characteristic polynomial of $B_3$ as
$$
\varphi_3(x)=x^{2}-(3\alpha s+s-\alpha-1)x+5\alpha s^{2}-2s^{2}-3\alpha s+s.
$$
Since the partition $V(G)=V(K_s)\cup V((2s-1)K_1)$ is equitable, it follows from Lemma 2.5 that $\rho_{\alpha}(G)$ is the largest root of
$\varphi_3(x)=0$. Note that $n=3s-1$. By a direct computation, we obtain
$$
\varphi_3(n-3)=\varphi_3(3s-4)=(4-4\alpha)s^{2}+(12\alpha-16)s-4\alpha+12.
$$

If $0\leq\alpha\leq\frac{1}{2}$, then $n=3s-1\geq f(\alpha)=14$. Thus, we infer $s\geq5$, and so
\begin{align*}
\varphi_3(n-3)=&(4-4\alpha)s^{2}+(12\alpha-16)s-4\alpha+12\\
\geq&5(4-4\alpha)s+(12\alpha-16)s-4\alpha+12\\
=&(4-8\alpha)s-4\alpha+12\\
\geq&5(4-8\alpha)-4\alpha+12\\
=&32-44\alpha\\
>&0.
\end{align*}

If $\frac{1}{2}<\alpha\leq\frac{2}{3}$, then $n=3s-1\geq f(\alpha)=17$. Thus, we have $s\geq6$, and so
\begin{align*}
\varphi_3(n-3)=&(4-4\alpha)s^{2}+(12\alpha-16)s-4\alpha+12\\
\geq&6(4-4\alpha)s+(12\alpha-16)s-4\alpha+12\\
=&(8-12\alpha)s-4\alpha+12\\
\geq&6(8-12\alpha)-4\alpha+12\\
=&60-76\alpha\\
>&0.
\end{align*}

If $\frac{2}{3}<\alpha\leq\frac{3}{4}$, then $n=3s-1\geq f(\alpha)=20$. Thus, we get $s\geq7$, and so
\begin{align*}
\varphi_3(n-3)=&(4-4\alpha)s^{2}+(12\alpha-16)s-4\alpha+12\\
\geq&7(4-4\alpha)s+(12\alpha-16)s-4\alpha+12\\
=&(12-16\alpha)s-4\alpha+12\\
\geq&7(12-16\alpha)-4\alpha+12\\
=&96-116\alpha\\
>&0.
\end{align*}

If $\frac{3}{4}<\alpha<1$, then $n=3s-1\geq f(\alpha)=\frac{5}{1-\alpha}+1$. Thus, we deduce $s\geq\frac{5}{3(1-\alpha)}+\frac{2}{3}$, and so
\begin{align*}
\varphi_3(n-3)=&(4-4\alpha)s^{2}+(12\alpha-16)s-4\alpha+12\\
\geq&(4-4\alpha)\left(\frac{5}{3(1-\alpha)}+\frac{2}{3}\right)s+(12\alpha-16)s-4\alpha+12\\
=&\frac{1}{3}((28\alpha-20)s-12\alpha+36)\\
\geq&\frac{1}{3}\left((28\alpha-20)\left(\frac{5}{3(1-\alpha)}+\frac{2}{3}\right)-12\alpha+36\right)\\
=&\frac{1}{9(1-\alpha)}(-20\alpha^{2}+92\alpha-32)\\
>&0.
\end{align*}

From the discussion above, we conclude $s\geq5$ and $\varphi_3(n-3)>0$. Next, we consider the derivative of $\varphi_3(x)$, we have
$\varphi_3'(x)=2x-(3\alpha s+s-\alpha-1)$. It follows from $s\geq5$ and $n=3s-1$ that
\begin{align*}
\varphi_3'(n-3)=&2(n-3)-(3\alpha s+s-\alpha-1)\\
=&(5-3\alpha)s+\alpha-7\\
\geq&5(5-3\alpha)+\alpha-7\\
=&18-14\alpha\\
>&0.
\end{align*}
Combining this with $\varphi_3(n-3)>0$ and $\eta(n)>n-3$, we conclude $\rho_{\alpha}(G)<n-3<\eta(n)$ for $s\geq5$, which contradicts
$\rho_{\alpha}(G)>\eta(n)$.

\noindent{\bf Subcase 2.2.2.} $i(G-S)=2s$.

In this subcase, $G=K_s\vee(K_{n_1}\cup2sK_1)$. If $n_1\geq2$, then $n=3s+n_1\geq3s+2$. By a similar proof as that of Subcase 2.1 above, we have
$\rho_{\alpha}(G)<\eta(n)$ for $2\leq s\leq\frac{n-2}{3}$, which contradicts $\rho_{\alpha}(G)>\eta(n)$.

If $n_1=0$, then $G=K_s\vee2sK_1$ and $n=3s$. Consider the partition $V(G)=V(K_s)\cup V(2sK_1)$. The corresponding quotient matrix of $A_{\alpha}(G)$
equals
\begin{align*}
B_4=\left(
  \begin{array}{ccc}
    2\alpha s+s-1 & 2s(1-\alpha)\\
    (1-\alpha)s &  \alpha s\\
  \end{array}
\right).
\end{align*}
The characteristic polynomial of $B_4$ is
$$
\varphi_4(x)=x^{2}-(3\alpha s+s-1)x+5\alpha s^{2}-2s^{2}-\alpha s.
$$
Since the partition $V(G)=V(K_s)\cup V(2sK_1)$ is equitable, it follows from Lemma 2.5 that $\rho_{\alpha}(G)$ is the largest root of
$\varphi_4(x)=0$. Notice that $n=3s$. By a simple calculation, we get
$$
\varphi_4(n-3)=\varphi_4(3s-3)=(4-4\alpha)s^{2}+(8\alpha-12)s+6.
$$

If $0\leq\alpha\leq\frac{2}{3}$, then $n=3s\geq f(\alpha)\geq14$. Thus, we have $s\geq5$, and so
\begin{align*}
\varphi_4(n-3)=&(4-4\alpha)s^{2}+(8\alpha-12)s+6\\
\geq&5(4-4\alpha)s+(8\alpha-12)s+6\\
=&(8-12\alpha)s+6\\
\geq&5(8-12\alpha)+6\\
=&46-60\alpha\\
>&0.
\end{align*}

If $\frac{2}{3}<\alpha\leq\frac{3}{4}$, then $n=3s\geq f(\alpha)=20$. Thus, we conclude $s\geq7$, and so
\begin{align*}
\varphi_4(n-3)=&(4-4\alpha)s^{2}+(8\alpha-12)s+6\\
\geq&7(4-4\alpha)s+(8\alpha-12)s+6\\
=&(16-20\alpha)s+6\\
\geq&7(16-20\alpha)+6\\
=&118-140\alpha\\
>&0.
\end{align*}

If $\frac{3}{4}<\alpha<1$, then $n=3s\geq f(\alpha)=\frac{5}{1-\alpha}+1$. Thus, we get $s\geq\frac{5}{3(1-\alpha)}+\frac{1}{3}$, and so
\begin{align*}
\varphi_4(n-3)=&(4-4\alpha)s^{2}+(8\alpha-12)s+6\\
\geq&(4-4\alpha)\left(\frac{5}{3(1-\alpha)}+\frac{1}{3}\right)s+(8\alpha-12)s+6\\
=&\frac{1}{3}((20\alpha-12)s+18)\\
\geq&\frac{1}{3}\left((20\alpha-12)\left(\frac{5}{3(1-\alpha)}+\frac{1}{3}\right)+18\right)\\
=&\frac{1}{9(1-\alpha)}(-20\alpha^{2}+78\alpha-18)\\
>&0.
\end{align*}

From the discussion above, we obtain $s\geq5$ and $\varphi_4(n-3)>0$. In what follows, we consider the derivative of $\varphi_4(x)$, we get
$\varphi_4'(x)=2x-(3\alpha s+s-1)$. According to $s\geq5$ and $n=3s$, we obtain
\begin{align*}
\varphi_4'(n-3)=&2(n-3)-(3\alpha s+s-1)\\
=&(5-3\alpha)s-5\\
\geq&5(5-3\alpha)-5\\
=&20-15\alpha\\
>&0.
\end{align*}
Combining this with $\varphi_4(n-3)>0$ and $\eta(n)>n-3$, we deduce $\rho_{\alpha}(G)<n-3<\eta(n)$ for $s\geq5$, which is a contradiction to
$\rho_{\alpha}(G)>\eta(n)$. This completes the proof of Theorem 1.2. \hfill $\Box$

\section{Extremal graphs}

In this section, we create a graph to claim that the bound on $A_{\alpha}$-spectral radius established in Theorem 1.1 is sharp.

\medskip

\noindent{\textbf{Theorem 4.1.}} Let $\alpha\in[0,1)$, and let $\eta(n)$ be the largest root of
$x^{3}-((\alpha+1)n+\alpha-4)x^{2}+(\alpha n^{2}+(\alpha^{2}-2\alpha-1)n-2\alpha+1)x-\alpha^{2}n^{2}+(5\alpha^{2}-3\alpha+2)n
-10\alpha^{2}+15\alpha-8=0$. Then $\rho_{\alpha}(K_{n-3}\vee K_1\vee \overline{K_2})=\eta(n)$ and $K_{n-3}\vee K_1\vee \overline{K_2}$ is not
a $P_{\geq2}$-factor covered graph if $n\geq f(\alpha)$, where
\[
f(\alpha)=\left\{
\begin{array}{ll}
14,&if \ \alpha\in[0,\frac{1}{2}];\\
17,&if \ \alpha\in(\frac{1}{2},\frac{2}{3}];\\
20,&if \ \alpha\in(\frac{2}{3},\frac{3}{4}];\\
\frac{5}{1-\alpha}+1,&if \ \alpha\in(\frac{3}{4},1).\\
\end{array}
\right.
\]

\medskip

\noindent{\it Proof.} We partition the vertex set of the graph $K_{n-3}\vee K_1\vee\overline{K_2}$ as $V(K_{n-3})\cup V(K_1)\cup V(\overline{K_2})$.
Then the quotient matrix of $A_{\alpha}(K_{n-3}\vee K_1\vee\overline{K_2})$ in terms of the partition $V(K_{n-3})\cup V(K_1)\cup V(\overline{K_2})$
can be written as
\begin{align*}
B(K_{n-3}\vee K_1\vee\overline{K_2})=\left(
  \begin{array}{ccc}
    \alpha(n-1) & (1-\alpha)(n-3) & 2(1-\alpha)\\
    1-\alpha & n+\alpha-4 & 0\\
    1-\alpha & 0 & \alpha\\
  \end{array}
\right),
\end{align*}
whose characteristic polynomial is $x^{3}-((\alpha+1)n+\alpha-4)x^{2}+(\alpha n^{2}+(\alpha^{2}-2\alpha-1)n-2\alpha+1)x-\alpha^{2}n^{2}+(5\alpha^{2}-3\alpha+2)n
-10\alpha^{2}+15\alpha-8$. Since the partition is equitable, it follows from Lemma 2.5 that the largest root $\eta(n)$ of $x^{3}-((\alpha+1)n+\alpha-4)x^{2}+(\alpha n^{2}+(\alpha^{2}-2\alpha-1)n-2\alpha+1)x-\alpha^{2}n^{2}+(5\alpha^{2}-3\alpha+2)n-10\alpha^{2}+15\alpha-8=0$ is equal to the $A_{\alpha}$-spectral
radius of the graph $K_{n-3}\vee K_1\vee\overline{K_2}$, that is to say, $\rho_{\alpha}(K_{n-3}\vee K_1\vee\overline{K_2})=\eta(n)$. Write $S=V(K_1)$.
Then $i(K_{n-3}\vee K_1\vee\overline{K_2}-S)=2>1=2|S|-1$. According to Lemma 2.2 (\romannumeral2), we know that the graph $K_{n-3}\vee K_1\vee\overline{K_2}$
is not a $P_{\geq2}$-factor covered graph. Theorem 4.1 is proved. \hfill $\Box$

\section*{Data availability statement}

My manuscript has no associated data.

\section*{Declaration of competing interest}

The authors declare that they have no conflicts of interest to this work.


\end{document}